\documentclass[11pt,twoside]{amsart}

\usepackage{amsfonts}
\usepackage{amsmath}
\usepackage{amssymb}
\usepackage{latexsym}

\title[Morita theory ]
{On Morita theory for self-dual modules}
\author{Wolfgang Willems}
\address{Wolfgang Willems,
\newline Institut f\"ur Algebra und Geometrie,
\newline Fakult\"at f\"ur Mathematik,
\newline Otto-von-Guericke-Universit\"at,
\newline 39016 Magdeburg,
\newline Germany}
\email{wolfgang.willems@ovgu.de}
\author{Alexander Zimmermann}
\address{Alexander Zimmermann
\newline Universit\'e de Picardie,
\newline Facult\'e de Math\'ematiques et LAMFA (UMR 6140 du CNRS),
\newline 33 rue St Leu,
\newline F-80039 Amiens Cedex 1,
\newline France}
\email{alexander.zimmermann@u-picardie.fr}
\date{May 13, 2007; revised March 25, 2008}

\newtheorem{Example}{Example}[section]
\newtheorem{Theorem1}[Example]{Theorem}

\newtheorem{Remark}[Example]{Remark}

\newtheorem{Corollary}[Example]{Corollary}

\newtheorem{Lemma}[Example]{Lemma}
\newtheorem{Definition}[Example]{Definition}
\newtheorem{Proposition}[Example]{Proposition}

\newcommand{\lra}{\longrightarrow}

\newcommand{\ra}{\rightarrow}
\newcommand{\sdp}{\times\kern-.2em\vrule height1.1ex depth-.05ex}
\newcommand{\epi}{\lra \kern-.8em\ra}

\newcommand{\dickebox}{{\vrule height5pt width5pt depth0pt}}

\newcommand{\rad}{\mbox{rad}}

 \normalbaselines
\thanks{The first author is supported by a German-French grant ``PROCOPE'' of the DAAD, the second by a French-German grant
``partenariat Hubert Curien PROCOPE dossier 14765WB''. Both authors
gratefully acknowledge  financial support.}

\begin{document}

\begin{abstract}
Let $G$ be a finite group and let $k$ be a field of characteristic
$p$. It is known that a $kG$-module $V$ carries a non-degenerate $G$-invariant
bilinear form $b$ if and only if $V$ is self-dual. We show that whenever
a Morita bimodule $M$ which induces an equivalence
between two blocks $B(kG)$ and $B(kH)$
of group algebras $kG$ and $kH$ is self-dual then
the correspondence preserves self-duality. Even more, if the bilinear form on
$M$ is symmetric then for $p$ odd the correspondence
preserves the geometric type of simple modules. In characteristic $2$ this holds also true
for projective modules.
\end{abstract}

\maketitle

\section{Introduction}

Let $k$ be a field and let $G$ be a finite group. It is well-known that
a $kG$-module $V$, simply called a $G$-module, carries a
non-degenerate $G$-invariant $k$-bilinear form if and only if
$V$ is isomorphic to its dual $V^* = Hom_k(V,k)$.
Up to an automorphism of $V$, a simple self-dual $G$-module
$V \cong V^*$ carries
exactly one $G$-invariant bilinear form different from the zero.
In case
$k$ is large enough, for instance $k$ algebraically  closed, this means that
such a form, if it exists, is unique up to a scalar.
If in addition the characteristic $p$ of $k$ is odd this form is symmetric or
antisymmetric. The case $p=2$ turns out to be more subtle. Note that in this case an
antisymmetric form is symmetric. By Fong's Lemma
(\cite{Huppert}, Chap. VII, Theorem 8.13),
a simple self-dual $G$-module $V$ different from the trivial module always
carries a non-degenerate $G$-invariant symplectic form. However, it may happen
that $V$ even has a non-degenerate $G$-invariant quadratic form. This holds
true, for instance, if $G$ is solvable but not in general.
What really happens is not an easy task to decide; it is a question
in cohomology (see \cite{Sin-Wi}), at least for arbitrary groups.
To be brief throughout
the paper we call a $G$-module of symmetric, antisymmetric, symplectic,
resp. quadratic
type if $V$ carries a non-degenerate $G$-invariant form of symmetric, antisymmetric,
symplectic, resp. quadratic type. The reader not familiar with duality theory of modules
may refer Chapter VII of \cite{Huppert}.

\medskip

In representation theory of finite groups $G$ we are often
interested in the category of modules
belonging to a  $p$-block $B_G$ of $kG$ where the underlying
field $k$ is of characteristic $p$.
Instead of investigating $B_G$ we may study a Morita equivalent $p$-block $B_H$
for some other group $H$.
Recall that one calls two blocks Morita equivalent if their module categories are
equivalent as exact categories. Morita theory shows that the module categories of the
two blocks are equivalent if and only if
there is an $H \times G$-bimodule
$M=\mbox{}_H M_G$ where $M$ is a progenerator for $B_H$ and for $B_G$.
Moreover the functor
$$ V \mapsto M \otimes_{kG} V $$
defines an equivalence from the module category of $B_G$
onto that of $B_H$.
The correspondence preserves all functorial properties of module categories,
in particular being simple, indecomposable,
projective, Loewy structure
and the like. For more details the reader is referred to  Chapter 4 in \cite{RMO}.

\medskip

In this note we study the behaviour of self-duality including the type of forms under
Morita equivalence. Given a group $\Gamma$ recall that a $\Gamma$-module $T$ carries
a non degenerate $\Gamma$-invariant bilinear form if and only if $T$ is self-dual.
A first step to our main result, useful in its own right, is the following.
\medskip

%\begin{Theorem}\label{iff}
\noindent
{\bf Theorem A.} {\it Let $B_G$ and $B_H$ be blocks for $G$ and $H$ respectively.
Suppose that the $H \times G$-bimodule $M=\mbox{}_HM_G$ defines a
Morita equivalence between $B_G$ and $B_H$ and suppose that $M$ is self-dual.
Then a $B_G$-left module $V$ is self-dual
if and only if $M\otimes_{kG}V$ is a self-dual $B_H$-left module. }
%\end{Theorem}

\medskip

We present now our main result for algebraically closed fields $k$ of odd characteristic.
Given a group $\Gamma$ we call a self-dual $k\Gamma$-module $T$ of
{\em symmetric type }
if $T$ carries a non degenerate $\Gamma$-invariant symmetric bilinear form.
If $T$ is simple and self-dual the $\Gamma$-invariant form is unique up
to a scalar and the form is then either symmetric
or antisymmetric. We call this the {\em type} of the simple module $T$.

\medskip

%\begin{Theorem} \label{Theo4.3}
\noindent
{\bf Theorem B.} {\it
Let $k$ be an algebraically closed field of odd characteristic.
Let $M$ be an $ H \times G$-bimodule inducing a Morita equivalence between
a $kG$-block $B_G$ and a $kH$-block $B_H$. If $M$ is of symmetric type then the
equivalence preserves the
type of simple self-dual modules. }
%\end{Theorem}

\medskip

The main result for fields of characteristic $2$ is weaker. We call a
$k\Gamma$-module $T$ of {\em quadratic type} if $T$ admits a $\Gamma$-invariant non degenerate quadratic form.

\medskip

%\begin{Theorem} \label{even2}
\noindent
{\bf Theorem C.} {\it
Let $k$ be an algebraically closed field of characteristic $2$.
Let $M$ be an $ H \times G$-bimodule inducing a Morita equivalence between
a $kG$-block $B_G$ and a $kH$-block $B_H$. If $M$ is of symmetric type
then the equivalence preserves
the property "being of quadratic type"
on projective  modules. }
%\end{Theorem}

\medskip

Theorem A, Theorem B and Theorem C
give invariants under Morita equivalence in case the
Morita bimodule which induces the equivalence is self-dual for Theorem A,
or even more of symmetric type for Theorem B and Theorem C. Moreover, we show that this is actually a property of the module
categories, rather than of the particular bimodule. Indeed Remark~\ref{remarkpicgroup}
shows that if $M$ is a self-dual Morita bimodule inducing an equivalence between
two blocks of group algebras, then any Morita bimodule between these blocks is self-dual.
If the Morita bimodule is moreover of symmetric type, then every Morita bimodule
between these blocks is self-dual.

Conditions on the Morita equivalence are actually necessary. The assertions of the theorems are false if we do not
require suitable conditions on the Morita module.
If the two blocks $B_G$ and $B_H$ of the group $G$ and the group $H$ are
Morita equivalent, it happens that
every simple module over $B_G$ is self-dual whereas no simple module over $B_H$ is
self-dual. In particular, no self-dual Morita bimodule exists in this case. An example
is given by a non principal $7$-block of the triple cover of the McLaughlin group
and the principal $7$-block of the Mathieu group $M_{23}$. For details see
Example \ref{Example2} below.
Even worse, two Morita equivalent blocks $B_G$ and $B_H$ may admit a self-dual Morita bimodule,
but the types of the forms on the simple modules of $B_G$ and of $B_H$ are different. Hence
it may happen that no Morita bimodule can be of symmetric type. An example
is given by a non principal $3$-block of $SL_2(5)$ and the principal $3$-block of
${\mathfrak A}_5$. Details are found in Example \ref{Example1} below.

The paper is organised as follows. In Section~\ref{Examplesection} the above-cited examples are
developed in detail. Section~\ref{selfdualbimodule} is the technical heart of the paper since there Proposition~\ref{selfdualMorita} is shown:
Given a self-dual, right-projective $H\times G$-bimodule $M$ and a self-dual $G$ module $V$,
Proposition~\ref{selfdualMorita} gives an explicit isomorphism
between the $H$-module $M\otimes_G V$ and its dual. The proposition applies as well for stable equivalences of Morita type. Except unicity, most of what follows is valid
in this more general setting. Moreover, Theorem A is proven there. Section~\ref{oddmorita} then proves Theorem B and Section~\ref{evenmorita} proves Theorem C.

\section{Examples}
\label{Examplesection}

In this section we show that Morita equivalence
does not need to preserve self-duality of simple modules. In case
it does,
examples show that the geometric type of simple modules
may vary via the correspondence.

In order to state the examples
we need that  blocks   are determined
by their Brauer trees up to Morita equivalence.
Kupisch determined the structure of indecomposable
modules of a Brauer tree algebra, and hence he implicitly determined their
module category, i.e., the Morita equivalence class of a Brauer tree algebra.
In modern terms Kupisch's result reads as follows.

\begin{Theorem1} {\rm (Kupisch \cite{Kupisch1,Kupisch2}) }
Let $A_1$ and $A_2$ be two Brauer tree algebras.
Then $A_1$ and $A_2$ are Morita equivalent if and only if their Brauer trees,
 including the exceptional vertex and its multiplicity, coincide.
\end{Theorem1}

In the following examples we shall use Brauer trees
computed by Hiss and Lux in \cite{HissLux}. The geometric type
of simple self-dual modules, i.e. symmetric or alternating, is taken from \cite{Conway} or \cite{Jansen}.

\begin{Example}\label{Example1} {\rm
Let $k$ be an algebraically closed field of characteristic $3$ and let
$G= SL_2(5)$ be the double cover of the alternating group ${\mathfrak A}_5$
of degree $5$.
Then a non principal $3$-block, say $B_1(G)$, of $G$ is a Brauer tree
algebra  corresponding to the
Brauer tree
\unitlength1cm
\begin{center}
\begin{picture}(4,1)
\put(0,.5){$T:$}
\put(.95,.5){$\bullet$}
\put(1.2,.6){\line(1,0){1}}
\put(2.3,.5){$\bullet$}
\put(2.55,.6){\line(1,0){1}}
\put(3.7,.5){$\bullet$}
\end{picture}
\end{center}
without an exceptional vertex.
The ordinary characters of the algebra have dimensions $2$, $4$ and $2$. Both
of the simple modules in characteristic $3$ of $B_1(G)$ admit an, up to scalar unique,
$G$-invariant bilinear
form. This form is alternating for both modules.

The principal $3$-block of the alternating group ${H=\mathfrak A}_5$ is also a Brauer tree algebra
with Brauer tree $T$.
The ordinary characters of this algebra have dimensions $1$, $5$ and $4$.
In characteristic $3$ both simple modules
admit an, up to scalar unique, ${\mathfrak A}_5$-invariant bilinear form which is
symmetric in both cases.}
\end{Example}

As a consequence, the type of bilinear forms is usually not preserved
by Morita equivalence, even if  self-duality
is preserved on simple modules.

\begin{Example}\label{Example2} {\rm
Let $k$ be again an algebraically closed field but of characteristic $7$ and let
$G=3McL$ be the triple cover of a the McLaughlin simple group $McL$.

Then, a non principal $7$-block, say $B_2(3McL)$, of $G$ is a Brauer tree
algebra corresponding to the Brauer tree
\unitlength1cm
\begin{center}
\begin{picture}(6,1)
\put(0,.5){$S:$}
\put(.95,.5){$\bullet$}
\put(1.2,.6){\line(1,0){1}}
\put(2.3,.5){$\bullet$}
\put(2.55,.6){\line(1,0){1}}
\put(3.6,.5){$\bullet$}
\put(3.85,.6){\line(1,0){1}}
\put(4.95,.5){$\circ_2$}
\end{picture}
\end{center}
with an exceptional vertex $\circ_2$ of multiplicity $2$.
The ordinary characters of the algebra have dimensions $792$, $4752$, $6336$
and $2376$ read from the left in the order of the tree.
None
of the three simple modules in characteristic $7$ is self-dual, hence none
of them admit a $G$-invariant bilinear form.

The principal $7$-block $B_3(M_{23})$ of the Mathieu group
$H=M_{23}$ is a Brauer tree algebra also with Brauer tree $S$. The
ordinary characters of the algebra have dimensions $1$, $1035$,
$2024$ and $990$  read from the left. The exceptional vertex
corresponds to the two ordinary characters of degree $990$. All
simple modules in characteristic $7$ admit an, up to scalar
unique, $H$-invariant bilinear form which is symmetric in all
three cases. }
\end{Example}

As a consequence every Morita bimodule
inducing a Morita equivalence between
$B_3(M_{23})$ and $B_2(3McL)$ has the property that
each self-dual simple
$B_3(M_{23})$-module is sent to a simple $B_2(3McL)$-module
which is not self-dual.

\section{Self-dual Morita bimodules}
\label{selfdualbimodule}

Let $H$ and $G$ be finite groups and let $k$ be a field.
Furthermore, let $M$ be an $H \times G$-bimodule where $H$ acts on the
left and $G$ on the right. Then $M^* = Hom_k(M,k)$ becomes an
$H \times G$-bimodule by setting
$$  (hfg)m = f(h^{-1}mg^{-1}) $$
for $h \in H, g \in G, f \in M^*$ and $m \in M$. Indeed,
the module structure follows by
$$ \begin{array}{rcl}
  ((h_2h_1)f(g_1g_2))m & = & f((h_2h_1)^{-1}m(g_1g_2)^{-1})\\[1ex]
            & = & f(h_1^{-1}h_2^{-1}mg_2^{-1}g_1^{-1})\\[1ex]
            & = & (h_1fg_1)(h_2^{-1}mg_2^{-1})\\[1ex]
            & = & (h_2(h_1fg_1)g_2)m.
  \end{array}
$$

\begin{Definition} Let $M$ be an $H \times G$-bimodule. A $k$-bilinear form
$(\cdot,\cdot)$ on $M$ is called $H \times G$-invariant if
$$ (hmg,hng) =(m,n)$$
for all $h \in H$, all $g \in G$ and all $m,n \in M$.
\end{Definition}

The next lemma is well-known. For the readers convenience we give the proof.

\begin{Lemma} \label{dualform}
 Let $M$ be an $H \times G$-bimodule. Then $M$ is isomorphic
to $M^*$ as an $H \times G$-bi\-mo\-dule if and only if $M$ carries a
non-degenerate $H \times G$-invariant bilinear form. More precisely,
if $\varphi:M\rightarrow M^*$ is an isomorphism, then $(m,m'):=\varphi(m)(m')$
is a non-degenerate $H \times G$-invariant bilinear form. If $(\cdot,\cdot)$
is a non-degenerate $H \times G$-invariant bilinear form, then the map
$M\ni m\mapsto (\cdot,m)\in M^*$ is an isomorphism of $H \times G$-bi\-mo\-dules.
\end{Lemma}

{\em Proof.} If $\varphi:M \rightarrow M^*$
denotes an $H \times G$-isomorphism we put
$$ (m,n) := \varphi(m)n$$
for $m,n \in M$. Clearly, the form $(\cdot,\cdot)$ is $k$-linear,
but it is also $H \times G$-invariant
since
$$ \begin{array}{rcl}
 (hmg,hng) & = & \varphi(hmg)(hng) \\[1ex]
           & = & (h\varphi(m)g)(hng) \qquad \quad (\mbox{since $\varphi$ is $H \times G$-linear})\\[1ex]
           & = & \varphi(m)(h^{-1}(hng)g^{-1}) \quad (\mbox{the module structure of $M^*$})\\[1ex]
           & = & \varphi(m)n \\[1ex]
           & = & (m,n).
  \end{array}
$$
In addition $(\cdot, \cdot)$ is non-degenerate since $\varphi$ is an isomorphism.

Conversely, if $(\cdot, \cdot)$ is a non-degenerate $H \times G$-invariant bilinear form on $M$
then we put
$$ \varphi(m)n := (m,n)$$
for all $m,n \in M$. Note that $\varphi$ is a $k$-isomorphism since
the bilinear form is non-degenerate. \\
Furthermore, $\varphi$ is $H \times G$-linear since
$$ (hmg,n) = \varphi(hmg)n$$
and
$$ (hmg,n)=(m,h^{-1}ng^{-1}) = \varphi(m)(h^{-1}ng^{-1}) = (h \varphi(m)g)n$$
for all $h \in H, g \in G$ and $m,n \in M$. \hfill\dickebox

\begin{Proposition} \label{selfdualMorita} Let $M$ be a self-dual $H \times G$-bimodule
which is projective as a $G$-module. If $V$ is a
self-dual $G$-left module then $M \otimes_{kG} V$ is a self-dual $H$-left
module.
More precisely, if $\beta:M\rightarrow M^*$ is an $H\times G$-linear isomorphism and
$\alpha:V\rightarrow V^*$ is a $G$-linear isomorphism, then
\begin{eqnarray*}
M \otimes_{kG} V&\longrightarrow&\left(M \otimes_{kG} V\right)^*\\
m\otimes v&\mapsto&\left(m'\otimes v'\mapsto
\sum_{g\in G}\beta(mg^{-1})(m')\cdot\alpha(v)(g^{-1}v')\right)
\end{eqnarray*}
is an $H$-linear isomorphism.
\end{Proposition}

{\em Proof.} Let $\beta: M \rightarrow M^*$ denote an $H \times G$-isomorphism
and $\alpha: V \rightarrow V^*$ a $G$-isomorphism. Therefore,
since  $\alpha$ and $\beta$ are $kG$-linear the map
$\gamma:=\beta\otimes_{kG}\alpha:M\otimes_{kG}V\rightarrow M^*\otimes_{kG}V^*$
is a well-defined $k$-vector space isomorphism satisfying
$(\beta\otimes_{kG}\alpha)(m\otimes_{kG} v)=\beta(m)\otimes_{kG}\alpha(v)$
for all $m\in M$ and $v\in V$ (cf. e.g. \cite[Theorem 12.10]{CRold}).

Now, let $m \in M, v \in V$ and $h \in H$. Since
$$ \begin{array}{rcl}
  \gamma(h(m \otimes_{kG} v)) & = & \gamma(hm \otimes_{kG} v) \\[1ex]
             & = & \beta(hm) \otimes_{kG} \alpha(v)  \\[1ex]
             & = & h\beta(m)\otimes_{kG}\alpha(v) \quad \mbox{(since $\beta$ is $H$-linear)} \\[1ex]
             & = & h\gamma(m\otimes_{kG}v).
  \end{array}
$$
the map  $\gamma$ is $H$-linear.

\medskip

Finally, we claim that
     $$ M^* \otimes_{kG} V^* \cong (M \otimes_{kG} V)^* $$
as $H$-left modules. We prove this via a sequence of isomorphisms.

\medskip

Observe that $V^*$ is not only a $G$-left but in a natural way
also a $G$-right module via
the definition
$$ (fg)v = f(gv)$$ for
$f \in V^*, g \in G$ and $v \in V$.

\medskip

Considering $V^*$ as $G$-right module $Hom_{kG}(M,V^*)$ makes sense.
With this observation we have

\medskip

{\bf Claim (i)}
$Hom_{kG}(M,V^*) \cong Hom_{k}(M \otimes_{kG} V, k) = (M \otimes_{kG} V)^* $ as $H$-left modules.\\[1ex]
The isomorphism is given by sending $\alpha \in Hom_{kG}(M,V^*)$
to the map $ \Phi({\alpha}): M \otimes_{kG} V
 \rightarrow k$ with $ \Phi({\alpha})(m \otimes_{kG} v) =\alpha(m)v$ where
 $m \in M$ and $v \in V$.
Note, that $ \Phi$ is well-defined since
$ \alpha$ is $G$-linear, for
$$\Phi(\alpha)(mg\otimes_{kG}v)=\alpha(mg)(v)=\left(\alpha(m)g\right)(v)=\alpha(m)(gv)=
\Phi(\alpha)(m\otimes_{kG}gv).$$
Furthermore, $\Phi$ is $H$-linear since
$$ \Phi({h\alpha})(m \otimes_{kG} v) = (h\alpha)(m)v = \alpha(h^{-1}m)v =
 \Phi({\alpha})(h^{-1}m \otimes_{kG} v) = (h \Phi({\alpha}))(m \otimes_{kG} v).$$
Finally, $\alpha \rightarrow \Phi({\alpha})$ is obviously a monomorphism.
In order to show that $\Phi$ is an isomorphism we define
a $k$-linear map
$$ \Psi: Hom_{k}(M \otimes_{kG} V, k) {\rightarrow} Hom_{kG}(M,V^*)$$
by putting $(\Psi(\varphi)(m))(v):=\varphi(m\otimes v)$ for  $\varphi\in
Hom_{k}(M \otimes_{kG} V, k)$,  $m\in M$ and $v\in V$.
Since $\varphi$ is defined over $M\otimes_{kG}V$
we have
$$ (\Psi(\varphi)(mg))v= \varphi(mg \otimes_{kG} v) = \varphi(m \otimes_{kG} gv) =
   (\Psi(\varphi)(m))(gv) = (\Psi(\varphi)(m)g)v $$
for $ m \in M, g \in G$ and $v \in V$.
Thus
$\Psi(\varphi)\in Hom_{kG}(M,V^*)$. Moreover,
$$((\Phi\circ \Psi)(\varphi))(m\otimes -)=\Phi(m\mapsto \varphi(m\otimes-))= \varphi(m\otimes-)  $$
for $m\in M$. Thus
$\Phi$ is surjective, hence an isomorphism.

\medskip

{\bf Claim (ii)} $Hom_{kG}(M,kG) \otimes_{kG} V^* \cong Hom_{kG}(M,V^*)$
as $H$-left modules:

\medskip

The isomorphism is given by
$$ \alpha \otimes_{kG} f \mapsto \widehat{\alpha \otimes_{kG} f} $$
for $\alpha \in Hom_{kG}(M,kG)$ and $f \in V^*$ where
$$ \left((\widehat{\alpha \otimes_{kG} f})(m)\right)(v) =
(f \alpha(m))(v)=f\left(\alpha(m)v\right)$$
for $ m \in M$ and $v\in V$. Here $V^*$ is considered as a
$G$-right module  by the natural action
as defined above.
Furthermore, the space $Hom_{kG}(M,kG)$ is an $H \times G$-bimodule via the action
$$(hfg)(m):=g^{-1} f(h^{-1}m)$$
for $m\in M$, $g\in G$, $h \in H$ and $f\in Hom_{kG}(M,kG)$.

\medskip

First we show that the map $\widehat{\phantom{xxx}}$ is well-defined.

\medskip
\noindent
Indeed, for $f\in V^*$, $g\in G$, $\alpha\in Hom_{kG}(M,kG)$, $v\in V$ and $m\in M$ one have
\begin{eqnarray*}
\left((\widehat{(\alpha g) \otimes_{kG} f})(m)\right)(v) &=&
\left(f ((\alpha g)(m))\right)(v)\ \qquad \quad \mbox{(definition of $\widehat{\phantom{xx}}$})\\
&=& f\left(\left((\alpha g)(m)\right)v\right) \ \ \qquad \quad \mbox{(right module structure of $V^*$)}\\
&=& f\left(\left(g^{-1}(\alpha (m))\right)v\right)\;\; \mbox{(right module structure of $Hom_{kG}(M,kG)$)}\\
&=& f\left(g^{-1}\left( \alpha(m) v\right)\right)\;\; \qquad \ \ \mbox{(associativity of the $kG$-action on $V^*$)}\\
&=&(gf)(\alpha(m)v)\;\; \qquad \quad  \quad \ \mbox{(left module structure of $V^*$})\\
&=&\left((\widehat{\alpha  \otimes_{kG} (gf)})(m)\right)(v)\;\;\mbox{(definition
of $\widehat{\phantom{xx}}$)}
\end{eqnarray*}

\medskip

Furthermore, the map $\widehat{\phantom{xx}}$ is
$H$-linear, for
$$ (\widehat{h\alpha \otimes_{kG} f})(m) =  f(h\alpha)(m) =
f  \alpha(h^{-1}m)  = (\widehat{\alpha \otimes_{kG} f })(h^{-1}m) =
(h(\widehat{\alpha \otimes_{kG} f}))(m).$$
Thus it remains to prove that $\widehat{\phantom{xx}}$ is a $k$-isomorphism. First we consider the case
$M=kG$. Let $ \alpha \in Hom_{kG}(kG,V^*)$, hence
$$ \alpha(\sum_{g \in G}a_gg) = \sum_{g \in G} a_g \alpha(g) =
\sum_{g \in G} a_g \alpha(1)g = f \sum_{g \in G} a_g g $$
where $\alpha(1) = f \in V^*$ and $a_g \in k$ for $g \in G$.
If $id$ denotes the identity in $Hom_{kG}(kG,kG)$ then
$$ (\widehat{id \otimes_{kG} f})(\sum_{g \in G} k_gg) = f\sum_{g \in G} k_gg
= \alpha(\sum_{g \in G} k_gg),$$
hence $\widehat{id \otimes_{kG} f} = \alpha$. Thus the above map is an epimorphism,
and therefore an isomorphism since the dimensions of the spaces in (ii) are obviously  equal for $M=kG$.
Since any free $G$-module is a direct sum of modules isomorphic $kG$ the
map $\widehat{\phantom{xx}}$ is also an isomorphism if $M$ is a free module.
Since a projective module is
a direct summand of a free module and since the isomorphisms are all
compatible with taking direct summands the proof of (ii) is complete.

\medskip

{\bf Claim (iii)} $Hom_k(M,k) \cong Hom_{kG}(M,kG)$ as
$H \times G$-bimodules: \\[1ex]
Abstractly this follows by  Frobenius reciprocity
\cite[Chapter VI, formula (8.7)]{MacLane} and the fact
that $kG$ is a symmetric algebra. However,
we need the isomorphism explicit. Let $(\cdot,\cdot)$
be the symmetrising bilinear form making $kG$ into a symmetric algebra.
For all $x,y\in kG$ the value of $(x,y)$ is the
coefficient of $1$ in $xy\in kG$ (cf. e.g. \cite[Theorem 62.1]{CRold}).
Then the map
\begin{eqnarray*}
Hom_{kG}(M,kG)&\stackrel{\rho_M}{\lra}&Hom_k(M,k)\\
f&\mapsto&(m\mapsto (f(m),1))
\end{eqnarray*}
is a $k$-linear.
If $\rho_M(f)=0$ we get
 $$ 0 = (f(mx),1)=(f(m)x,1)=(f(m),x)$$
for all $ m \in M$ and all $x \in kG$.
Since $(\cdot,\cdot)$ is non degenerate
$f(m)=0$ for all $m\in M$ and hence $f=0$. Therefore $\rho_M$ is injective.

Moreover, $\rho_M$ is a morphism of $H \times G$-bimodules.
Indeed, for  $g\in G$ and $h\in H$, for $f\in Hom_{kG}(M,kG)$ and
$m\in M$, we have
\begin{eqnarray*}
\left(\rho_M(hfg)\right)(m)&=&((hfg)(m),1)\,\, \qquad \mbox{ (definition of $\rho_M$)}\\
&=&(g^{-1}f(h^{-1}m),1) \quad \mbox{ ($H \times G$-bimodule structure of $Hom_{kG}(M,kG)$)}\\
&=&(1,g^{-1}f(h^{-1}m))\,\, \ \ \mbox{ (since $(\cdot,\cdot)$ is symmetric)}\\
&=&(g^{-1},f(h^{-1}m))\,\,\quad \,\mbox{ (since $(\cdot,\cdot)$ is associative)}\\
&=&(f(h^{-1}m),g^{-1})\,\,\ \ \ \,\mbox{ (since $(\cdot,\cdot)$ is symmetric)}\\
&=&(f(h^{-1}m)g^{-1},1)\,\, \ \ \mbox{ (since $(\cdot,\cdot)$ is associative)}\\
&=& (f(h^{-1}mg^{-1}),1)\,\, \ \ \ \mbox{ ($f$ is $kG$-linear)}\\
&=&\rho_M(f)(h^{-1}mg^{-1})\,\,\,\mbox{ (definition of $\rho_M$)}\\
&=&\left(h\left(\rho_M(f)\right)g\right)(m)\,\,\,
\mbox{ ($H \times G$-bimodule structure of $Hom_{k}(M,k)$)}
\end{eqnarray*}
and so $\rho_M(hfg)=h\left(\rho_M(f)\right)g$ which shows that $\rho_M$ is $H \times G$-linear.

\medskip

The inverse isomorphism is given by  sending $ \varphi \in
Hom_k(M,k)$ to
 $\Lambda_M({\varphi}) \in Hom_{kG}(M,kG)$ defined by
$$ \Lambda_M({\varphi})(m) = \sum_{x \in G} \varphi(mx^{-1})x$$
for $m\in M$:
First note, that $\Lambda_M(\varphi)$ is $G$-linear since it is the trace
under the $G$-action.  Furthermore,
\begin{eqnarray*}
\left((\rho_M\circ\Lambda_M)(\varphi)\right)(m)&=&(\sum_{x \in G} \varphi(mx^{-1})x,1)\\
&=&\sum_{x \in G}\varphi(mx^{-1})(x,1)\\
&=&\sum_{x \in G}\varphi(mx^{-1}) \delta_{x,1}\\
&=&\varphi(m)
\end{eqnarray*}
for $m \in M$ and $\varphi \in Hom_k(M,k)$
where $\delta$ denotes the Kronecker symbol. Thus, $\rho_M$ is
surjective and $\rho_M^{-1}=\Lambda_M$ is an isomorphism
of $H \times G$-bimodules.

\bigskip

Since $\Lambda_M$ is $H\times G$-linear, $\Lambda_M\otimes_{kG}id_{V^*}$
is well defined (cf. e.g. \cite[Theorem 12.10]{CRold}).
Summarising (i), (ii) and (iii) we get explicit isomorphisms
$$ \begin{array}{rcl}
 M\otimes_{kG}V & \stackrel{\beta\otimes\alpha}{\longrightarrow} & M^* \otimes_{kG} V^* \\[1ex]
                & \stackrel{\Lambda_M\otimes id_{V^*}}{\longrightarrow} &  Hom_{kG}(M,kG) \otimes_{kG} V^* \\[1ex]
                & \stackrel{\widehat{\phantom{xxx}}}{\longrightarrow} & Hom_{kG}(M,V^*) \\[1ex]
                & \stackrel{\Phi}{\longrightarrow} & Hom_{k}(M \otimes_{kG} V,k) =
                (M \otimes_{kG} V)^*
    \end{array}
$$
as $H$-left modules.

We need to verify the formula given in the Proposition.
\begin{eqnarray*}
m\otimes v&\stackrel{\beta\otimes\alpha}{\mapsto}&\beta(m)\otimes_{kG}\alpha(v)\\
&\stackrel{\Lambda_M\otimes id}{\mapsto}&\sum_{g\in G}\beta(mg^{-1})g\otimes_{kG}\alpha(v)\\
& =&\sum_{g\in G}\beta(mg^{-1})\otimes_{kG} g\alpha(v)\\
&\stackrel{\widehat{\phantom{xxx}}}{\mapsto}&
m'\mapsto
\left(
v'\mapsto
\left(
\sum_{g\in G} g\alpha(v)\left(\beta(mg^{-1})(m') v'\right)
\right)
\right)\\
&=&
m'\mapsto
\left(
v'\mapsto
\left(
\sum_{g\in G}\beta(mg^{-1})(m')\cdot g\alpha(v)(v')
\right)
\right)\\
& \stackrel{\Phi}{\mapsto} &\left((m'\otimes v')\mapsto
\left(\sum_{g\in G}\beta(mg^{-1})(m')\cdot\alpha(v)(g^{-1} v')\right)\right)
\end{eqnarray*}
This complets the proof. \hfill\dickebox

\bigskip

The following is the translation of the isomorphism given in Proposition~\ref{selfdualMorita}
in terms of bilinear forms.

\begin{Corollary} \label{corollary1} Let $M$ be an $H \times G$-bimodule
which is projective
as a $G$-module. Suppose that $M$ carries a non-degenerate
$H \times G$-invariant non-degenerate
bilinear form $b$. Let $V$ be a $G$-left module with a non-degenerate $G$-invariant
bilinear form $B$. Then, for  $m, m' \in M$ and $v, v' \in V$,
$$  \tilde{B}(m \otimes_{kG} v, m' \otimes_{kG} v') = \sum_{g \in G} b(mg^{-1},m')B(v,g^{-1}v')$$
 defines a  non-degenerate
$H$-invariant bilinear form on $M \otimes_{kG} V$.
\end{Corollary}

{\em Proof.} By Lemma \ref{dualform}, the bilinear form $b$ is associated to
an $H \times G$-isomorphism \\$\beta: M \rightarrow M^*$ via
$$ b(m,m') = \beta(m)m'$$
for $m,m' \in M$. Similarly, $B$ is associated to a $G$-isomorphism
$\alpha:V \rightarrow V^*$
via
$$ B(v,v') = \alpha(v)v'$$
for $v, v' \in V$.
Proposition~\ref{selfdualMorita} shows that
an isomorphism
$$\Psi: M\otimes_{kG}V \longrightarrow \left(M\otimes_{kG}V\right)^*$$
is given by
$$ \begin{array}{rcl}
\left(\Psi(m\otimes v)\right)(m'\otimes v')&=&
 \sum_{g \in G} \beta(mg^{-1})(m') \cdot\alpha(v) (g^{-1}v') \\[1ex]
 & = & \sum_{g \in G} b(mg^{-1},m')B(v,g^{-1}v').
\end{array}
 $$
Thus, by Lemma \ref{dualform}, the  bilinear form $\tilde{B}$
defined on $M \otimes_{kG} V$ by
$$\tilde{B}(m \otimes_{kG} v, m' \otimes_{kG} v') =
\sum_{g \in G} b(mg^{-1},m') B(v,g^{-1}v')$$
is non-degenerate and $H$-invariant.
\hfill\dickebox

\medskip

{\em Proof of Theorem A.} By Proposition~\ref{selfdualMorita}, the functor $M \otimes_{kG} -$
maps self-dual $G$-modules to self-dual $H$-modules.

On the other hand, $Hom_{kH}(M,kH)$ becomes a $G \times H$-bimodule
via
$$ (g\alpha h)(m) = \alpha(mg)h$$
for all $g \in G,h \in H$ and $\alpha \in Hom_{kH}(M,kH)$. Then $Hom_{kH}(M,kH)$
is the inverse Morita bimodule to $M$ and hence $Hom_{kH}(M,kH)$ defines
the inverse functor from the category of $B_H$-modules onto
the category of $B_G$-modules. In order to prove the only if part it is,
by Proposition~\ref{selfdualMorita}, enough to show
that the  $G \times H$-bimodule $Hom_{kH}(M,kH)$ is self-dual as well.

First note that $$Hom_{kH}(M,kH) \cong Hom_k(M,k)= \mbox{}_G(M^*)_H$$ as
 $G \times H$-bimodules (by Claim (iii) in the
 proof of Proposition~\ref{selfdualMorita},
 changing left- and right structures there).
If we define a $G \times H$-structure on $M$ by setting
$$   g \circ m \circ h = h^{-1}mg^{-1}$$
then the inverse of the given isomorphism $ \mbox{}_HM_G \longrightarrow \mbox{}_HM^*_G$
defines
an isomorphism from
$\mbox{}_G(M^*)_H$ to $M \cong (\mbox{}_G(M^*)_H)^*$
and we are done. \hfill\dickebox

\bigskip

Example~\ref{Example2} shows that not all
Morita bimodules do admit a non-degenerate invariant bilinear form.

\begin{Remark}\label{remarkpicgroup} {\rm
Given two Morita equivalences between two blocks $B_G$ and $B_H$
of group algebras $kG$ and $kH$,
the composition of the one with the inverse of the other gives a
Morita self-equivalence of the block $B_G$.
Now, the isomorphism classes of Morita self-equivalences of $B_G$ form
the Picard group $Pic(B_G)$ which is isomorphic to the outer automorphism
group $Out(B_G):=Aut_{k-\mbox{\scriptsize alg}}(B_G)/Inn(B_G)$ where
as usual $Inn(B_G)$ denotes the inner automorphisms of $B_G$. Each
$ \alpha \in Aut_{k-\mbox{\scriptsize alg}}(B_G)$ induces an equivalence
by sending a module $V$ to $\mbox{}^\alpha V$, which is defined to be $V$ as a
$k$-vector space, but
where the action of $B_G$ on $\mbox{}^\alpha V$ is given by
$$ a \bullet v = \alpha(a) v$$
for $ a \in B_G$ and $v \in V$.
Since $(\mbox{}^\alpha V)^* \cong \mbox{}^\alpha (V^*)$
the automorphism $\alpha$ preserves the property of self-duality on modules.
Since $B_H\simeq End_{B_G}(M)$, whatever isomorphism is taken to identify
$B_H$ with $End_{B_G}(M)$ the question whether a given module $M\otimes_{kG}V$
carries a non degenerate $H$-invariant bilinear
form is independent of this choice. }
\end{Remark}

There are Morita equivalences which are induced by  self-dual bimodules.

\begin{Example} \label{example2} {\rm
Let $G$ be a finite $p$-solvable group and let $N=O_{p'}(G)$.
We put $H=G/N$ and $$e = \frac{1}{|O_{p'}(G)|}\sum_{g\in O_{p'}(G)}g.$$
Clearly, $e$ is the block idempotent of the principal block $B_0$ of $kG$
and $ekG \cong kG/N = kH$ as $H \times G$-bimodules where $N$ obviously acts
trivially on both sides. The bimodule $ \mbox{}_HM_G =ekG$ is a Morita bimodule
which induces an equivalence between the principal block of $kG$ and the
group algebra $kH$. In fact
$$End_{kG}(ekG) \cong ekGe \cong kG/N =kH,$$ and obviously $ekG$
is a progenerator of $ekG =B_0$. Moreover, $ekG \cong (ekG)^*$ as $H \times G$-bimodules
since, by a well-known result on duality,
$$(ekG)^* \cong \bar{e}kG = ekG$$
where $\mbox{}^{-}$ denotes the anti-automorphism of $kG$ defined by
$ g \rightarrow \bar{g} = g^{-1}$.
Thus the equivalence preserves self-duality, by
Proposition~\ref{selfdualMorita}.
Clearly, this fact is obvious since $N$ acts
trivially on all modules in $B_0$.
 }
\end{Example}

\section{Morita equivalence and forms in odd characteristic}
\label{oddmorita}

Throughout this section let $k$ be a field of odd characteristic
(which includes the characteristic 0 case as well).
A $G$-module $V$ is called of {\it symmetric} resp. {\it antisymmetric type}
if $V$ carries a non-degenerate $G$-invariant symmetric resp. antisymmetric
bilinear form.

\begin{Proposition}\label{forms}
Let $M$ be an $ H \times G$-bimodule inducing a Morita equivalence between
a $kG$-block $B_G$ and a $kH$-block $B_H$. Suppose that
$M$ is of symmetric type.
Let $V$ be a self-dual $G$-left module.
If $V$ is of symmetric resp. antisymmetric type then $M \otimes_{kG} V$ is
of symmetric resp. antisymmetric type.
\end{Proposition}

{\em Proof.} Suppose $b$ is a
symmetric bilinear form and the bilinear form $B$ is symmetric resp. antisymmetric
if $V$ is of symmetric resp. antisymmetric type.
With the notation of Corollary \ref{corollary1}, we have
$$ \begin{array}{rcl}
  \tilde{B}(m \otimes_{kG} v, m' \otimes_{kG} v') & = &
  \sum_{g \in G} b(mg^{-1},m')B(v,g^{-1}v') \\[1ex]
                  & = & \sum_{g \in G} b(m,m'g)B(gv,v') \qquad
                  \mbox{(since $b$ and $B$ are $G$-invariant)} \\[1ex]
                   & = & \pm \sum_{g \in G} b(m'g,m)B(v',gv)
                   \quad \mbox{(by symmetry of $b$ and $B$)} \\[1ex]
                   & = &  \pm \tilde{B}(m' \otimes_{kG} v', m \otimes_{kG} v),
   \end{array}
$$
where the $+$ sign appears for symmetric and the $-$ sign
 for antisymmetric type modules $V$. \mbox{} \qquad \qquad
\hfill\dickebox

\medskip

{\em Proof of Theorem B.} This is an immediate consequence of Proposition \ref{forms}. \hfill\dickebox

\begin{Example} \label{example3} {\rm We consider again Example \ref{example2}. Note that the
 Morita module $ekG$ is
of symmetric type. As a form one can take the standard non-degenerate
$G$-invariant symmetric form on $kG$ defined by
$$ b(g,g') = \delta_{g,g'}$$
for $g,g' \in G$ restricted to $ekG$. The restriction is indeed
non-degenerate as one can check easily. Thus Proposition
\ref{forms} applies.}
\end{Example}

Now let $V$ be a simple self-dual $kG$-module.
Thus $V$ carries a $G$-invariant non-degenerate
bilinear form. If $k$ is algebraically closed then - up to scalars - there is
exactly one non-degenerate $G$-invariant form which is symmetric or
antisymmetric since the
characteristic is odd (\cite{Huppert},Chap. VII, 8.12). Thus, for
$k$ algebraically closed, a simple self-dual $G$-module is either of
symmetric or antisymmetric type but not of both.
In particular, the Morita module $M=ekG$ in Example \ref{example3} is never of
antisymmetric type since otherwise, by Proposition \ref{forms}, there are
simple $H$-modules which are of symmetric and antisymmetric type.
Indeed it is well-known, that a group algebra $kH \cong ekG$ in odd
characteristic never has a non-degenerate
$H$-invariant antisymmetric form.

\begin{Remark} {\rm Let $k$ be an algebraically closed field of odd characteristic.
By definition, a Morita bimodule $M=\mbox{}_HM_G$ is a projective $kG$-right and
a projective $kH$-left module. If $M$ carries a non-degenerate $H \times G$-invariant
symmetric form
then a projective indecomposable $kG$- resp. $kH$-module of antisymmetric type occurs
with even
multiplicity in $M$. Note that, by (\cite{Wil}, Proposition 2.2), an indecomposable projective module is of antisymmetric
type if and only if its head is of antisymmetric type.}
\end{Remark}

\section{Morita equivalence and forms in even characteristic}
\label{evenmorita}

In characteristic $2$ the geometry of self-dual modules turns out to be more subtle.
In this case an alternating form is symmetric.
By Fong's Lemma (\cite{Huppert}, Chap. VII, Theorem 8.13), a simple
non-trivial self-dual $kG$-module $V$ always carries a non-degenerate $G$-invariant symplectic
form, but it may happen that $V$ even has a non-degenerate $G$-invariant quadratic form.
This happens if $H^1(G,V)=0$ (\cite{Sin-Wi}, Proposition 2.4) and
a connecting homomorphism to the first cohomology group
carries the essential information what really happens.
We call an arbitrary module $V \cong V^*$ of quadratic type if it carries a
non-degenerate $G$-invariant quadratic form, say $Q$.
The corresponding $G$-invariant non-degenerate symplectic form $B$ defined by
$$ B(v,v') = Q(v+v') - Q(v) -Q(v')$$
for $v, v' \in V$ is called the Polarisation of $Q$. Note that the quadratic form $Q$
is non-degenerate if
$$  \rad \, Q  := \{ v \mid v \in \rad \, B, \  Q(v)=0  \} =0. $$
Throughout the following $k$ is assumed to be a perfect field of characteristic $2$.
The next Lemma can be found as Proposition 2.2 in \cite{Gow-Wil}.

\begin{Lemma} \label{even1} Let $P$ be a projective $kG$-module not containing
the projective cover of the trivial module. Then each $G$-invariant symmetric
form on $P$ is the Polarisation of a $G$-invariant quadratic form on $P$.
\end{Lemma}

We are now ready to prove Theorem C.

\medskip

{\em Proof of Theorem C.} Let $P$ be a projective $kG$-module of quadratic type. Then
$P$ is an orthogonal sum of orthogonal indecomposable modules. By
(\cite{Wil}, Lemma 3.6), an orthogonal indecomposable module is
either indecomposable or of the form $S \oplus S^*$ where $S$ is
indecomposable. Since $$ (M \otimes_{kG} (S \oplus S^*) )^* \cong
M^* \otimes_{kG} (S \oplus S^*)^* \cong M \otimes_{kG}
         (S \oplus S^*)$$
by Claim (i),  (ii) and  (iii) in the proof of
Proposition~\ref{selfdualMorita} we need to prove the Theorem only for orthogonal indecomposable projective
modules.

So let $P$ be an orthogonal indecomposable projective $kG$-module in $B_G$ of quadratic type.
Thus $P$ carries a non-degenerate $G$-invariant symplectic form $B$ which is also symmetric.
According to Corollary \ref{corollary1} the form $\tilde{B}$ defined by
$$  \tilde{B}(m \otimes_{kG} p, m' \otimes_{kG} p') = \sum_{g \in G} b(mg^{-1},m') B(p,g^{-1}p')$$
for $m, m' \in M$ and $p, p' \in P$ is a non-degenerate $H$-invariant bilinear form on
$M \otimes_{kG} P$. Since $b$ and $B$ are symmetric the form
$\tilde{B}$ is symmetric as well. Thus, if the projective  module
$M \otimes_{kG} P$ does not contain the
projective cover of the trivial $kH$-module, say $P_H(1)$, then $M \otimes_{kG} P$ is of quadratic type, by Lemma \ref{even1}.
In case it contains $P_H(1)$ it is either $P_H(1)$ or a direct sum of two copies of $P_H(1)$.
But, by
(\cite{Wil}, Remark 3.5 (b)), $P_H(1)$ and hence any finite direct sum of it is always of quadratic type if $2 \mid |H|$.
In case $2 \nmid |H|$ the trivial module $k$
is projective and obviously of quadratic type. This completes the proof.
\hfill\dickebox

\bigskip

At the moment we are not able to prove an analogue of Theorem C for simple
self-dual modules with the methods we have developed so far.

\end{document}